\documentclass{article}
\usepackage{amsmath,amsfonts,amssymb,bm,hyperref,verbatim}
\numberwithin{equation}{section}
\newtheorem{theorem}{Theorem}[section]
\newtheorem{lemma}{Lemma}[section]

\newtheorem{remark}{Remark}[section]

\usepackage[T1]{fontenc}
\usepackage[utf8]{inputenc}
\usepackage{authblk}

\author{Olga Chervova\thanks{O.Chervova@ucl.ac.uk, \url{http://www.homepages.ucl.ac.uk/\~ucahoch/}} }
\author{Robert J.~Downes\thanks{R.Downes@ucl.ac.uk, \url{http://www.homepages.ucl.ac.uk/\~zcahc37/}} }
\author{Dmitri Vassiliev\thanks{D.Vassiliev@ucl.ac.uk, \url{http://www.homepages.ucl.ac.uk/\~ucahdva/}}}
\affil{Department of Mathematics,
University College London,\\ Gower Street, London WC1E 6BT, UK}

\begin{document}

\title{Spectral theoretic characterization of the massless Dirac operator}
\maketitle
\begin{abstract}
We consider an elliptic self-adjoint first order differential
operator acting on pairs (2-columns) of complex-valued half-densities over
a connected compact 3-dimensional manifold without boundary.
The principal symbol of our operator is assumed to be trace-free.
We study the spectral function which is the sum of squares of Euclidean norms
of eigenfunctions evaluated at a given point of the manifold,
with summation carried out over all eigenvalues between zero
and a positive~$\lambda$. We derive an explicit two-term asymptotic formula
for the spectral function as $\lambda\to+\infty$, expressing the
second asymptotic coefficient via the trace of the subprincipal symbol and the
geometric objects encoded within the principal symbol --- metric, torsion
of the tele\-parallel connection and topological charge.
We then address the question: is our operator a massless Dirac operator
on half-densities?
We prove that
it is a massless Dirac operator on half-densities if and only if the following two
conditions are satisfied at every point of the manifold: a)~the
subprincipal symbol is proportional to the identity matrix and
b)~the second asymptotic coefficient of the spectral function is zero.
\end{abstract}

\tableofcontents

\section{Main results}
\label{Main results}

Consider a first order differential
operator $A$ acting on 2-columns
$v=\begin{pmatrix}v_1&v_2\end{pmatrix}^T$
of complex-valued half-densities
over a connected compact 3-dimensional manifold $M$ without boundary.
(See subsection 1.1.5 in \cite{mybook} for definition of half-density.)
We assume the coefficients of the operator $A$ to be infinitely smooth. We also
assume that the operator $A$ is formally self-adjoint (symmetric):
\begin{equation}
\label{formally self-adjoint}
\int_Mw^*Av\,dx=\int_M(Aw)^*v\,dx
\end{equation}
for all infinitely smooth
$v,w:M\to\mathbb{C}^2$. Here and further on
the superscript~$\,{}^*\,$ in matrices, rows and columns
indicates Hermitian conjugation in $\mathbb{C}^2$
and $dx:=dx^1dx^2dx^3$, where $x=(x^1,x^2,x^3)$ are local
coordinates on $M$.

Let $A_1(x,\xi)$ be the principal symbol of the operator $A$,
i.e.~matrix obtained by leaving in $A$ only the leading (first order)
derivatives and replacing each $\partial/\partial x^\alpha$ by
$i\xi_\alpha$, $\alpha=1,2,3$.
Here $\xi=(\xi_1,\xi_2,\xi_3)$ is the variable dual to the position
variable $x$; in physics literature the $\xi$ would be referred to
as \emph{momentum}.
Our principal symbol $A_1(x,\xi)$ is a $2\times 2$
Hermitian matrix-function on the cotangent bundle $T^*M$, linear in
every fibre $T_x^*M$ (i.e.~linear in $\xi$).

Throughout this paper we assume that the principal symbol
$A_1(x,\xi)$ is trace-free for all $(x,\xi)\in T^*M$
and that
\begin{equation}
\label{ellipticity condition}
\det A_1(x,\xi)\ne0,\qquad\forall(x,\xi)\in T'M,
\end{equation}
where $T'M:=T^*M\setminus\{\xi=0\}$
(cotangent bundle with the zero section removed).
The assumption (\ref{ellipticity condition}) is a version of the ellipticity condition.

Under the above assumptions $A$ is a self-adjoint operator in
$L^2(M;\mathbb{C}^2)$ (Hilbert space of square integrable
complex-valued column ``functions'') with domain $H^1(M;\mathbb{C}^2)$
(Sobolev space of complex-valued column ``functions'' which are
square integrable together with their first partial derivatives) and
the spectrum of $A$ is discrete, with eigenvalues accumulating
to $\pm\infty$.
Let $\lambda_k$ and
$v_k=\begin{pmatrix}v_{k1}(x)&v_{k2}(x)\end{pmatrix}^T$ be
the eigenvalues and eigenfunctions of the operator $A$. The
eigenvalues $\lambda_k$ are enumerated in increasing order with
account of multiplicity,
using a positive index $k=1,2,\ldots$ for positive $\lambda_k$
and a nonpositive index $k=0,-1,-2,\ldots$ for nonpositive $\lambda_k$.

We will be studying the \emph{spectral function} and the \emph{counting function}.
The spectral function is the real density defined as
\begin{equation}
\label{definition of spectral function}
e(\lambda,x,x):=\sum_{0<\lambda_k<\lambda}\|v_k(x)\|^2,
\end{equation}
where $\|v_k(x)\|^2:=[v_k(x)]^*v_k(x)$ is the square of the
Euclidean norm of the eigenfunction $v_k$ evaluated at the point
$x\in M$ and $\lambda$ is a positive parameter (spectral parameter).
The counting function is the function
\begin{equation}
\label{definition of counting function}
N(\lambda):=\,\sum_{0<\lambda_k<\lambda}1\ =\int_Me(\lambda,x,x)\,dx\,.
\end{equation}
In other words, $N(\lambda)$ is the number of eigenvalues $\lambda_k$
between zero and $\lambda$.

We aim to derive, under appropriate assumptions on Hamiltonian
trajectories, two-term asymptotics for the spectral function
(\ref{definition of spectral function})
and the counting function
(\ref{definition of counting function}),
i.e.~formulae of the type
\begin{equation}
\label{two-term asymptotic formula for spectral function}
e(\lambda,x,x)=a(x)\,\lambda^3+b(x)\,\lambda^2+o(\lambda^2),
\end{equation}
\begin{equation}
\label{two-term asymptotic formula for counting function}
N(\lambda)=a\lambda^3+b\lambda^2+o(\lambda^2)
\end{equation}
as $\lambda\to+\infty$, where the real constants $a$, $b$ and real densities
$a(x)$, $b(x)$ are related in accordance with
\begin{equation}
\label{a via a(x)}
a=\int_Ma(x)\,dx,
\end{equation}
\begin{equation}
\label{b via b(x)}
b=\int_Mb(x)\,dx.
\end{equation}

In our recent paper \cite{jst_part_a} we performed a comprehensive analysis of two-term spectral asymptotics
for general first order elliptic systems. In doing this we showed that
all previous publications on systems gave
formulae for the second asymptotic coefficient that were either
incorrect or incomplete (i.e.~an algorithm for the calculation
of the second asymptotic coefficient rather than an explicit formula),
see Section 11 of \cite{jst_part_a} for the appropriate bibliographic review.
The correct formula for the coefficient $b(x)$ was the main result of \cite{jst_part_a}.

The problem examined in the current paper is a special case of that from~\cite{jst_part_a}.
Namely, in the current paper we make the following additional
assumptions as compared to \cite{jst_part_a}:
\begin{equation}
\label{assumption manifold has dimension 3}
\text{our manifold has dimension 3},
\end{equation}
\begin{equation}
\label{assumption number of equations in our system is 2}
\text{the number of equations in our system is 2},
\end{equation}
\begin{equation}
\label{assumption operator is differential}
\text{our operator is differential (as opposed to pseudodifferential)},
\end{equation}
\begin{equation}
\label{assumption principal symbol is trace-free}
\text{the principal symbol is trace-free}.
\end{equation}
The need for a detailed analysis of the special case
(\ref{assumption manifold has dimension 3})--(\ref{assumption principal symbol is trace-free})
is driven by applications to the massless Dirac operator.

The additional assumptions
(\ref{assumption manifold has dimension 3})--(\ref{assumption principal symbol is trace-free})
lead to the following simplifications as compared to \cite{jst_part_a}.
\begin{itemize}
\item
The subprincipal symbol $A_\mathrm{sub}$ does not depend
on the dual variable $\xi$ (momentum) and is a function of $x$
(position) only. Recall that the subprincipal symbol is the zeroth order term
of the full symbol of the first order operator $A$ written in a way
which makes it  invariant under coordinate transformations,
see formula (\ref{definition of subprincipal symbol}) for formal
definition and subsection 2.1.3 in \cite{mybook}
for background material.
\item
The principal symbol $A_1$ admits a geometric description.
\end{itemize}
The first of these simplifications is trivial whereas the second is not. We list below
the geometric objects encoded within the principal symbol.

\

\textbf{Geometric object 1: the \emph{metric.}}
Observe that the determinant of the principal symbol is a negative definite quadratic form
\begin{equation}
\label{definition of metric}
\det A_1(x,\xi)=-g^{\alpha\beta}\xi_\alpha\xi_\beta
\end{equation}
and the coefficients $g^{\alpha\beta}(x)=g^{\beta\alpha}(x)$,
$\alpha,\beta=1,2,3$, appearing in (\ref{definition of metric})
can be interpreted as
components of a (contravariant) Riemannian metric.
This implies, in
particular, that our Hamiltonian (positive eigenvalue of the
principal symbol) takes the form
\begin{equation}
\label{Hamiltonian expressed via metric}
h^+(x,\xi)=\sqrt{g^{\alpha\beta}(x)\,\xi_\alpha\xi_\beta}
\end{equation}
and the $x$-components of our Hamiltonian trajectories become
geodesics.

\

\textbf{Geometric object 2: the \emph{teleparallel connection.}}
This is an affine connection
defined as follows. Suppose we have a covector $\xi$ based at
the point $x\in M$ and we want to construct a parallel covector
$\tilde\xi$ based at the point $\tilde x\in M$. This is done by solving the
linear system of equations
\begin{equation}
\label{definition of parallel transport}
A_1(\tilde x,\tilde\xi)=A_1(x,\xi).
\end{equation}
Equation (\ref{definition of parallel transport}) is equivalent to a
system of three real linear algebraic equations for the three
real unknowns, components of the covector $\tilde\xi$, and it is easy to
see that this system has a unique solution. It is also easy to see that
the affine connection defined by formula (\ref{definition of parallel transport})
preserves the Riemannian norm of covectors,
i.e.~$g^{\alpha\beta}(\tilde x)\,\tilde\xi_\alpha\tilde\xi_\beta
=g^{\alpha\beta}(x)\,\xi_\alpha\xi_\beta$, hence, it is metric
compatible. The parallel transport defined by formula
(\ref{definition of parallel transport}) does not depend on the
curve along which we transport the (co)vector, so our connection has
zero curvature. The word ``teleparallel'' (parallel at a distance)
is used in theoretical physics \cite{cartantorsionreview} to describe metric compatible affine
connections with zero curvature.
The origins of this terminology go back to the works of A.~Einstein and
\'E.~Cartan \cite{unzicker-2005-,MR2276051,MR543192}, though Cartan preferred
to use the term ``absolute parallelism'' rather than ``teleparallelism''.

The teleparallel connection coefficients
$\Gamma^\alpha{}_{\beta\gamma}(x)$ can be written down explicitly in
terms of the principal symbol,
see formula (\ref{formula for teleparallel connection coefficients}),
and this allows us to define yet another geometric object --- the torsion tensor
\begin{equation}
\label{definition of torsion}
T^\alpha{}_{\beta\gamma}
:=\Gamma^\alpha{}_{\beta\gamma}-\Gamma^\alpha{}_{\gamma\beta}\,.
\end{equation}
Further on we raise and lower indices of the torsion tensor using
the metric.

\

\textbf{Geometric object 3: the \emph{topological charge.}}
It turns out, see Section~\ref{Teleparallel connection}, that
the existence of a principal symbol implies that our manifold $M$ is parallelizable.
Parallelizability implies orientability.
Having chosen a particular orientation,
we allow only changes of local coordinates $x^\alpha$, $\alpha=1,2,3$, which preserve
orientation.

We define the topological charge as
\begin{equation}
\label{definition of relative orientation}
\mathbf{c}:=-\frac i2\sqrt{\det g_{\alpha\beta}}\,\operatorname{tr}
\bigl((A_1)_{\xi_1}(A_1)_{\xi_2}(A_1)_{\xi_3}\bigr),
\end{equation}
with the subscripts $\xi_\alpha$
indicating partial derivatives.
We show in Section \ref{Teleparallel connection} that
the number $\mathbf{c}$ defined by formula
(\ref{definition of relative orientation})
can take only two values, $+1$ or $-1$,
and describes the orientation of the principal symbol
relative to the chosen orientation of local coordinates.

\

We have identified three geometric objects encoded within the
principal symbol --- metric, teleparallel connection and topological charge.
Consequently, one would expect the coefficient $b(x)$
from formula
(\ref{two-term asymptotic formula for spectral function})
to be expressed via these three geometric objects and the subprincipal
symbol. This assertion is confirmed by the following theorem.

\begin{theorem}
\label{theorem 1.1}
The coefficients in the two-term asymptotics
(\ref{two-term asymptotic formula for spectral function})
are given by the formulae
\begin{equation}
\label{formula for a(x)}
a(x)=\frac1{6\pi^2}
\,\sqrt{\det g_{\alpha\beta}(x)}\,,
\end{equation}
\begin{equation}
\label{formula for b(x)}
b(x)=\frac1{8\pi^2}
\bigl(
[\,
3\,\mathbf{c}*\!T^\mathrm{ax}-2\operatorname{tr}A_\mathrm{sub}
\,]
\,\sqrt{\det
g_{\alpha\beta}}
\,\bigr)(x)\,,
\end{equation}
where
\begin{equation}
\label{definition of axial torsion}
T^\mathrm{ax}_{\alpha\beta\gamma}:=
\frac13(T_{\alpha\beta\gamma}+T_{\gamma\alpha\beta}+T_{\beta\gamma\alpha})
\end{equation}
is \emph{axial torsion}
(totally antisymmetric piece of the torsion tensor)
and $*$ is the Hodge star (\ref{definition of Hodge star}).
\end{theorem}

\begin{remark}
\label{remark on mollification}
The spectral and counting functions admit two-term asymptotic expansions
(\ref{two-term asymptotic formula for spectral function})
and
(\ref{two-term asymptotic formula for counting function})
only under appropriate assumptions on geodesic loops and closed geodesics respectively,
see Theorems 8.3 and 8.4 in \cite{jst_part_a}.
However, one can easily reformulate
asymptotic formulae
(\ref{two-term asymptotic formula for spectral function})
and
(\ref{two-term asymptotic formula for counting function})
in such a way that they remain valid without assumptions on geodesics:
this can easily be achieved, say, by taking a convolution with a function from Schwartz space
$\mathcal{S}(\mathbb{R})$, see Theorems 7.1 and 7.2 in \cite{jst_part_a}.
Thus, the second asymptotic coefficients of the spectral and counting functions are well-defined irrespective
of how many geodesic loops or closed geodesics we have. We introduced the second
asymptotic coefficients $b(x)$ and $b$ via the unmollified asymptotic expansions
(\ref{two-term asymptotic formula for spectral function})
and
(\ref{two-term asymptotic formula for counting function})
simply for the sake of clarity of presentation.
\end{remark}

The proof of Theorem \ref{theorem 1.1} is given in Sections
\ref{Reduction from the general setting}--\ref{Integration}.

\

We now turn our attention to the massless Dirac operator. This
operator is defined in Appendix \ref{The massless Dirac operator},
see formula (\ref{definition of Weyl operator}),
and it does not fit into our scheme because it is an operator
acting on a 2-component complex-valued spinor (Weyl spinor) rather than
a pair of complex-valued half-densities. However, on a parallelizable manifold
components of a spinor can be identified with half-densities. We
call the resulting operator
\emph{the massless Dirac operator on half-densities}. The explicit
formula for the massless Dirac operator on half-densities is
(\ref{definition of Weyl operator on half-densities}).

The massless Dirac operator on half-densities is an operator of the type
we are considering in this paper, i.e.~a self-adjoint first order elliptic
differential operator acting on 2-columns of complex-valued half-densities
and with a trace-free principal symbol. We address the
question: is a given operator $A$ a massless Dirac operator? The
answer is given by the following theorem which
is our main result.

\begin{theorem}
\label{theorem 1.2}
The operator $A$ is a massless Dirac operator on half-densities if and
only if the following two conditions are satisfied at every point of the
manifold $M$: a) the subprincipal symbol of the operator,
$A_\mathrm{sub}(x)$, is proportional to the identity matrix and b)
the second asymptotic coefficient of the spectral function, $b(x)$,
is zero.
\end{theorem}

Note that conditions a) and b) in Theorem \ref{theorem 1.2}
are invariant under special unitary transformations,
i.e.~transformations of the operator
\begin{equation}
\label{unitary transformation of operator A}
A\mapsto RAR^*,
\end{equation}
where
$R:M\to\mathrm{SU}(2)$
is an arbitrary smooth special unitary matrix-function.
The invariance of condition b) is obvious.
In fact, condition b)
is invariant under the action of a broader group:
the unitary matrix-function $R(x)$ appearing in formula
(\ref{unitary transformation of operator A})
does not have to be special.
As to condition a), its invariance is established by
examination of formula (9.3) from \cite{jst_part_a}
with the use of the special commutation properties of
trace-free Hermitian $2\times2$ matrices
(the anticommutator of a pair of trace-free Hermitian $2\times2$ matrices
is a multiple of the identity matrix).
The fact that the conditions of Theorem \ref{theorem 1.2}
are $\mathrm{SU}(2)$ invariant
is not surprising as the massless Dirac operator is designed
around the concept of $\mathrm{SU}(2)$ invariance,
see Property 4 in Appendix \ref{The massless Dirac operator}.

The proof of Theorem \ref{theorem 1.2} is given in Sections
\ref{The subprincipal symbol} and \ref{Proof of Theorem}.

\

Theorems \ref{theorem 1.1} and \ref{theorem 1.2} tell us that for the massless
Dirac operator on half-densities formulae
(\ref{two-term asymptotic formula for spectral function})
and
(\ref{two-term asymptotic formula for counting function})
read
\begin{equation}
\label{two-term asymptotic formula for spectral function for dirac}
e(\lambda,x,x)=\frac{\sqrt{\det g_{\alpha\beta}(x)}}{6\pi^2}\,\lambda^3+o(\lambda^2),
\end{equation}
\begin{equation}
\label{two-term asymptotic formula for counting function for dirac}
N(\lambda)=\frac{\operatorname{Vol}M}{6\pi^2}\,\lambda^3+o(\lambda^2),
\end{equation}
where $\operatorname{Vol}M$ is the volume of the Riemannian 3-manifold $M$.

\begin{remark}
\label{remark on spectral function desnity versus scalar}
The factor $\sqrt{\det g_{\alpha\beta}(x)}$ appears in the RHS of
(\ref{two-term asymptotic formula for spectral function for dirac})
because we are working with the massless Dirac operator on half-densities
(\ref{definition of Weyl operator on half-densities}) rather than with the
massless Dirac operator on spinors
(\ref{definition of Weyl operator}).
For the massless Dirac operator on spinors the spectral function
is a scalar field (as opposed to a density) and formula
(\ref{two-term asymptotic formula for spectral function for dirac})
reads
$e(\lambda,x,x)=\frac{1}{6\pi^2}\lambda^3+o(\lambda^2)$.
\end{remark}

\section{Reduction from the general setting}
\label{Reduction from the general setting}

As explained in Section \ref{Main results},
the problem considered in the current paper is a special case of
that from \cite{jst_part_a}.
Formulae (1.23) and (1.24) from \cite{jst_part_a} in our case read
\begin{equation}
\label{formula for a(x) version 1}
a(x)=
\int\limits_{h^+(x,\xi)<1}{d{\hskip-1pt\bar{}}\hskip1pt}\xi\,,
\end{equation}
\begin{equation}
\label{formula for b(x) sum of two}
b(x)=
b_1(x)+b_2(x)\,,
\end{equation}
where
\begin{equation}
\label{formula for b_1(x) version 1}
b_1(x)=-3
\int\limits_{h^+(x,\xi)<1}
(
[v^+]^*A_\mathrm{sub}v^+
)(x,\xi)\,
{d{\hskip-1pt\bar{}}\hskip1pt}\xi\,,
\end{equation}
\begin{equation}
\label{formula for b_2(x) version 1}
b_2(x)=\frac{3i}2
\int\limits_{h^+(x,\xi)<1}
\{
[v^+]^*,A_1-2h^+I,v^+
\}
(x,\xi)\,
{d{\hskip-1pt\bar{}}\hskip1pt}\xi\,.
\end{equation}
Here $h^+(x,\xi)$ is the positive eigenvalue of the principal symbol (see
also formula (\ref{Hamiltonian expressed via metric})),
$v^+(x,\xi)$ is the corresponding normalized eigenvector (2-column),
${d{\hskip-1pt\bar{}}\hskip1pt}\xi$ is shorthand for
${d{\hskip-1pt\bar{}}\hskip1pt}\xi:=(2\pi)^{-3}\,d\xi
=(2\pi)^{-3}\,d\xi_1d\xi_2d\xi_3$
and $I$ is the $2\times2$ identity matrix.
Curly brackets in formula
(\ref{formula for b_2(x) version 1})
denote the Poisson bracket on matrix-functions
\begin{equation}
\label{Poisson bracket on matrix-functions}
\{P,R\}:=P_{x^\alpha}R_{\xi_\alpha}-P_{\xi_\alpha}R_{x^\alpha}
\end{equation}
and its further generalization
\begin{equation}
\label{generalised Poisson bracket on matrix-functions}
\{P,Q,R\}:=P_{x^\alpha}QR_{\xi_\alpha}-P_{\xi_\alpha}QR_{x^\alpha}\,,
\end{equation}
with the subscripts $x^\alpha$ and $\xi_\alpha$
indicating partial derivatives and
the repeated tensor index $\alpha$ indicating summation over $\alpha=1,2,3$.

Put $P^+(x,\xi):=[v^+(x,\xi)][v^+(x,\xi)]^*$, which is the orthogonal projection onto
the eigenspace $\,\operatorname{span}v^+\,$ of the principal symbol.
We have $A_1-2h^+I=2h^+P^+-3h^+I$ and $\{[v^+]^*,P^+,v^+\}=0$, so formula
(\ref{formula for b_2(x) version 1}) can be rewritten as
\begin{equation}
\label{formula for b_2(x) version 2}
b_2(x)=-\frac{9i}2
\int\limits_{h^+(x,\xi)<1}
(
h^+
\{
[v^+]^*,v^+
\}
)(x,\xi)\,
{d{\hskip-1pt\bar{}}\hskip1pt}\xi\,.
\end{equation}

Our aim now is to evaluate the integrals
(\ref{formula for a(x) version 1}),
(\ref{formula for b_1(x) version 1})
and
(\ref{formula for b_2(x) version 2})
explicitly.

Formulae (\ref{formula for a(x) version 1}) and (\ref{Hamiltonian expressed via metric})
immediately imply (\ref{formula for a(x)}).

In order to evaluate the integral
(\ref{formula for b_1(x) version 1})
we rewrite this formula as
\[
b_1(x)=-3
\int\limits_{h^+(x,\xi)<1}
\operatorname{tr}
(
A_\mathrm{sub}P^+
)(x,\xi)\,
{d{\hskip-1pt\bar{}}\hskip1pt}\xi
\]
and use the fact that $P^+(x,\xi)=\frac1{2h^+(x,\xi)}(A_1(x,\xi)+h^+(x,\xi)\,I)$. We get
\[
b_1(x)=-3
\int\limits_{h^+(x,\xi)<1}
\frac1{2h^+(x,\xi)}
\operatorname{tr}
(
A_\mathrm{sub}(A_1+h^+I)
)(x,\xi)\,
{d{\hskip-1pt\bar{}}\hskip1pt}\xi\,.
\]
But $A_\mathrm{sub}$ does not depend on $\xi$
whereas $A_1$ and $h^+$ are, respectively, odd and even in $\xi$,
so the term $\frac1{2h^+}\operatorname{tr}(A_\mathrm{sub}A_1)$
integrates to zero, leaving us with
\begin{equation}
\label{formula for b_1(x) version 4}
b_1(x)=-\frac32(\operatorname{tr}A_\mathrm{sub})(x)
\int\limits_{h^+(x,\xi)<1}
{d{\hskip-1pt\bar{}}\hskip1pt}\xi
\,=-\frac1{4\pi^2}
\bigl(\operatorname{tr}A_\mathrm{sub}\,\sqrt{\det g_{\alpha\beta}}\,\bigr)(x)\,.
\end{equation}

In order to complete the proof of Theorem \ref{theorem 1.1} we need to evaluate
explicitly the integral (\ref{formula for b_2(x) version 2}). The next three
sections deal with this nontrivial issue.

\section{Teleparallel connection}
\label{Teleparallel connection}

We show in this section that the principal symbol generates a
teleparallel connection which allows us to reformulate the results
of our spectral analysis in a much clearer geometric language.

Let us show first that
the existence of a principal symbol implies that our manifold $M$ is parallelizable.
The principal symbol $A_1(x,\xi)$ is linear in $\xi$ so it can be written as
\begin{equation}
\label{principal symbol via Pauli matrices}
A_1(x,\xi)=\sigma^\alpha(x)\,\xi_\alpha\,,
\end{equation}
where $\sigma^\alpha(x)$, $\alpha=1,2,3$, are some trace-free Hermitian
$2\times2$ matrix-functions.
Let us denote the elements of the matrices $\sigma^\alpha$ as $\sigma^\alpha{}_{\dot ab}\,$,
where the dotted index, running through the values $\dot1,\dot2$, enumerates the rows
and the undotted index, running through the values $1,2$, enumerates the columns;
this notation is taken from \cite{MR2670535}. Put
\begin{equation}
\label{frame via Pauli matrices}
e_1{}^\alpha(x):=\operatorname{Re}\sigma^\alpha{}_{\dot12}(x),
\quad
e_2{}^\alpha(x):=-\operatorname{Im}\sigma^\alpha{}_{\dot12}(x),
\quad
e_3{}^\alpha(x):=\sigma^\alpha{}_{\dot11}(x).
\end{equation}
Formula (\ref{frame via Pauli matrices}) defines a triple of smooth real vector fields
$e_j(x)$, $j=1,2,3$, on the manifold $M$.
These vector fields are linearly independent at every point $x$ of the manifold:
this follows from formula (\ref{ellipticity condition}).
Thus, the triple of vector fields $e_j$ is a \emph{frame}.
The existence of a frame means that the manifold $M$ is parallelizable.

Conversely, given a frame $e_j$
we uniquely recover the principal symbol
$A_1(x,\xi)$ via formulae
(\ref{principal symbol via Pauli matrices}),
(\ref{Pauli matrices 1})
and
(\ref{Pauli matrices 2}).
Thus, a principal symbol is equivalent to a frame.
Of course, this equivalence statement relies on our \emph{a priori}
assumptions
(\ref{formally self-adjoint}),
(\ref{ellipticity condition})
and
(\ref{assumption manifold has dimension 3})--(\ref{assumption principal symbol is trace-free}).

It is easy to see that the frame elements $e_j$ are
orthonormal with respect to the metric
(\ref{definition of metric}). Moreover, the
metric can be determined directly from the frame as
\begin{equation}
\label{definition of metric via frame}
g^{\alpha\beta}=
\delta^{jk}e_j{}^\alpha\,e_k{}^\beta\,,
\end{equation}
where the repeated frame indices $j$ and $k$ indicate summation over $j,k=1,2,3$.
The two definitions of the metric,
(\ref{definition of metric})
and
(\ref{definition of metric via frame}),
are equivalent.

Parallelizability implies orientability, see Proposition 13.5 in \cite{lee}.
Having chosen a particular orientation,
we allow only changes of local coordinates $x^\alpha$, $\alpha=1,2,3$, which preserve
orientation and
define the Hodge star
in the standard way: the action of $\,*\,$ on a rank
$q$ antisymmetric tensor $Q$ is
\begin{equation}
\label{definition of Hodge star}
(*Q)_{\gamma_{q+1}\ldots\gamma_3}:=(q!)^{-1}\,\sqrt{\det g_{\alpha\beta}}\,
Q^{\gamma_1\ldots\gamma_q}\varepsilon_{\gamma_1\ldots\gamma_3}\,,
\end{equation}
where $\varepsilon$ is the totally antisymmetric quantity,
$\varepsilon_{123}:=+1$, and $g$ is the Riemannian metric
(\ref{definition of metric}).
Here and further on we identify differential forms with covariant antisymmetric tensors.
We raise and lower tensor indices using our metric.

Substituting formulae
(\ref{principal symbol via Pauli matrices})
and
(\ref{frame via Pauli matrices})
into
(\ref{definition of relative orientation})
we get
\begin{equation}
\label{definition of relative orientation more natural}
\mathbf{c}=\operatorname{sgn}\det e_j{}^\alpha.
\end{equation}
Formula (\ref{definition of relative orientation more natural})
provides an equivalent (and more natural) definition of topological charge.
It also explains why the topological charge,
initially defined in Section~\ref{Main results} in accordance with formula
(\ref{definition of relative orientation}),
can only take values $+1$ or $-1$.

The concept of a teleparallel connection
was already defined in Section~\ref{Main results} in accordance with
formula (\ref{definition of parallel transport}).
This connection can be equivalently defined via the frame as follows.
Suppose we have a vector $v$ based at
the point $x\in M$ and we want to construct a parallel vector
$\tilde v$ based at the point $\tilde x\in M$. We decompose the vector $v$
with respect to the frame at the point $x$, $v=c^je_j(x)$,
and reassemble it with the same coefficients $c^j$ at the point $\tilde x$,
defining $\tilde v:=c^je_j(\tilde x)$.

We now define the covariant derivative corresponding to the teleparallel
connection. Our teleparallel connection is a special case of an affine
connection, so we are looking at a covariant derivative acting on
vector/covector fields in the usual manner
\[
\nabla_\mu v^\alpha=\partial v^\alpha/\partial x^\mu+\Gamma^\alpha{}_{\mu\beta}\,v^\beta\,,
\qquad
\nabla_\mu w_\beta=\partial w_\beta/\partial x^\mu-\Gamma^\alpha{}_{\mu\beta}\,w_\alpha\,.
\]
The teleparallel connection coefficients are defined from the conditions
\begin{equation}
\label{conditions for teleparallel connection coefficients}
\nabla_\mu e_j{}^\alpha=0\,,
\end{equation}
where the $e_j$ are elements of our frame.
Formula (\ref{conditions for teleparallel connection coefficients}) gives a system
of 27 linear algebraic equations for the determination of 27 unknown
connection coefficients. It is known
(see, for example, formula (A2) in \cite{MR2573111}),
that the unique solution of this system is
\begin{equation}
\label{formula for teleparallel connection coefficients}
\Gamma^\alpha{}_{\mu\beta}=e_k{}^\alpha(\partial e^k{}_\beta/\partial x^\mu)\,,
\end{equation}
where
\begin{equation}
\label{definition of coframe}
e^k{}_\beta:=\delta^{kj}g_{\beta\gamma}e_j{}^\gamma.
\end{equation}
The triple of covector fields $e^k$, $k=1,2,3$, is called the \emph{coframe}.
The frame and coframe uniquely determine each other via the relation
\begin{equation}
\label{alternative definition of coframe}
e_j{}^\alpha e^k{}_\alpha=\delta_j{}^k.
\end{equation}

Note that our notation for the frame and coframe is taken from \cite{solovej}.
We feel it necessary to mention this because there is a whole range of different notation
for frames/coframes in mathematics and theoretical physics literature, which makes
the subject somewhat confusing.

One can check by performing explicit calculations that the teleparallel connection
has the following two important properties:
\begin{equation}
\label{teleparallel connection property 1}
\nabla_\alpha g_{\beta\gamma}=0,
\end{equation}
which means that the connection is metric compatible, and
\begin{equation}
\label{teleparallel connection property 2}
(\nabla_\alpha\nabla_\beta-\nabla_\beta\nabla_\alpha)v^\gamma=0
\ \,\text{for any vector field}\ \,v\,,
\end{equation}
which means that the Riemann curvature tensor is zero.
Properties
(\ref{teleparallel connection property 1})
and
(\ref{teleparallel connection property 2})
are the defining properties of a teleparallel connection:
a teleparallel connection is, by definition \cite{cartantorsionreview},
an affine connection satisfying
(\ref{teleparallel connection property 1})
and
(\ref{teleparallel connection property 2}).

The tensor characterizing the ``strength'' of the teleparallel
connection is not the Riemann curvature tensor but the torsion
tensor (\ref{definition of torsion}).
The teleparallel connection is, in a sense, the opposite of the
more common Levi-Civita connection:
the Levi-Civita connection has zero torsion but nonzero curvature,
whereas the teleparallel connection has nonzero torsion but zero curvature.
In our paper we distinguish these two affine connections by using
different notation for connection coefficients:
we write the teleparallel connection coefficients
as $\Gamma^\alpha{}_{\beta\gamma}$
and the Levi-Civita connection coefficients (Christoffel symbols)
as $\left\{{{\alpha}\atop{\beta\gamma}}\right\}$,
see formula (\ref{definition of Christoffel symbols}).
It is known, see formula (7.34) in \cite{nakahara},
that the two sets of connection coefficients are related as
$
\Gamma^\alpha{}_{\beta\gamma}
=
\left\{{{\alpha}\atop{\beta\gamma}}\right\}
+\frac12
(T^\alpha{}_{\beta\gamma}+T_\beta{}^\alpha{}_\gamma+T_\gamma{}^\alpha{}_\beta)
$.

Substituting
(\ref{formula for teleparallel connection coefficients})
into
(\ref{definition of torsion})
we arrive at the following explicit formula for the torsion tensor
of the teleparallel connection
\begin{equation}
\label{explicit formula for torsion}
T=e_j\otimes de^j\,,
\end{equation}
where the $d$ stands for the exterior derivative.
For the sake of clarity we rewrite formula
(\ref{explicit formula for torsion}) in more detailed form, retaining all tensor indices,
\begin{equation}
\label{more explicit formula for torsion}
T^\alpha{}_{\beta\gamma}
=e_j{}^\alpha
(
\partial e^j{}_\gamma/\partial x^\beta
-
\partial e^j{}_\beta/\partial x^\gamma
)\,.
\end{equation}
As always, the repeated index $j$ appearing in formulae
(\ref{explicit formula for torsion})
and
(\ref{more explicit formula for torsion})
indicates summation over $j=1,2,3$.

Torsion is a rank three tensor antisymmetric in the last
two indices. Because we are working in dimension three, it is
convenient, as in \cite{rotational_elasticity},
to apply the Hodge star in the last two indices
and deal with the rank two tensor
\begin{equation}
\label{definition of torsion with a star}
\overset{*}T{}^\alpha{}_\beta:=
\frac12\,T^{\alpha\gamma\delta}\,\varepsilon_{\gamma\delta\beta}
\,\sqrt{\det g_{\mu\nu}}
\end{equation}
instead.
Substituting
(\ref{explicit formula for torsion})
into
(\ref{definition of torsion with a star})
we get
\begin{equation}
\label{explicit formula for torsion with a star}
\overset{*}T=e_j\otimes\operatorname{curl}e^j\,,
\end{equation}
where
\begin{equation}
\label{definition of curl}
(\operatorname{curl}e^j)_\beta:=(*de^j)_\beta
=\frac12\,(de^j)^{\gamma\delta}\,\varepsilon_{\gamma\delta\beta}
\,\sqrt{\det g_{\mu\nu}}\,.
\end{equation}

\section{Relation between curvature of the $\mathrm{U}(1)$ connection
and torsion of the teleparallel connection}
\label{Relation}

This section is devoted to the examination of the integrand in formula
(\ref{formula for b_2(x) version 2}). Recall that the
curly brackets in this integrand denote the Poisson bracket on matrix-functions
(\ref{Poisson bracket on matrix-functions}).

As explained in Section 5 of \cite{jst_part_a}, the expression
$-i\{[v^+]^*,v^+\}$ is the scalar curvature of the $\mathrm{U}(1)$ connection
generated by the eigenspace $\,\operatorname{span}v^+\,$ of the principal symbol.
This curvature term
appears in the general setting of a first order elliptic system.
A feature of the particular case
(\ref{assumption manifold has dimension 3})--(\ref{assumption principal symbol is trace-free})
considered in the current paper is that the
scalar curvature of the $\mathrm{U}(1)$ connection can be expressed via
torsion of the teleparallel connection. This is a substantial simplification.
The teleparallel connection is a simpler geometric object than
the $\mathrm{U}(1)$ connection because the coefficients of
the teleparallel connection do not depend on the dual variable
(momentum), i.e.~they are ``functions'' on the base manifold $M$.
The relationship between the two connections is established by the following lemma.

\begin{lemma}
\label{Relation lemma}
The scalar curvature of the $\mathrm{U}(1)$ connection is
expressed via the torsion of the teleparallel connection, metric and topological charge as
\begin{equation}
\label{curvature via torsion and metric}
-i\{[v^+]^*,v^+\}(x,\xi)
=\frac{\mathbf{c}}2\,
\frac
{\overset{*}T{}^{\alpha\beta}(x)\,\xi_\alpha\xi_\beta}
{(g^{\mu\nu}(x)\,\xi_\mu\xi_\nu)^{3/2}}\,.
\end{equation}
\end{lemma}

\

Recall that the topological charge
$\mathbf{c}=\pm1$ is defined
in accordance with formula (\ref{definition of relative orientation})
or, equivalently, in accordance with
formula (\ref{definition of relative orientation more natural}).

\

\emph{Proof of Lemma \ref{Relation lemma}\ }
We give the proof for the case
\begin{equation}
\label{topological invariant equals one}
\mathbf{c}=+1\,.
\end{equation}
There is no need to give a separate proof for
the case $\mathbf{c}=-1$ as the two cases reduce to one another
by means of a) the observation that torsion
(\ref{explicit formula for torsion})
is invariant under inversion of the frame
and b) the identity
\begin{equation}
\label{sum of curvatures is zero}
\{[v^+]^*,v^+\}+\{[v^-]^*,v^-\}=0,
\end{equation}
where $v^-(x,\xi)$ is the normalized eigenvector of the principal symbol
corresponding to the negative eigenvalue.
Formula (\ref{sum of curvatures is zero}) is a special case of
formula (1.22) from \cite{jst_part_a}.

We fix an arbitrary point $Q\in T'M$ and prove formula
(\ref{curvature via torsion and metric}) at this point.
As the LHS and RHS of (\ref{curvature via torsion and metric}) are invariant under
changes of local coordinates~$x$,
it is sufficient to prove formula
(\ref{curvature via torsion and metric}) in
Riemann normal coordinates, i.e.~local coordinates
such that $x=0$ corresponds to the projection of the point $Q$ onto the base manifold,
$g_{\mu\nu}(0)=\delta_{\mu\nu}$ and $\frac{\partial g_{\mu\nu}}{\partial x^\lambda}(0)=0$.
Moreover, as the formula we are proving involves only
first partial derivatives in $x$, we may assume, without loss of generality,
that
\begin{equation}
\label{special metric}
g_{\mu\nu}(x)=\delta_{\mu\nu}
\end{equation}
for all $x$ in some neighbourhood of the origin.
In other words, it is sufficient to prove formula
(\ref{curvature via torsion and metric})
for the case of Euclidean metric.

As both the LHS and RHS of (\ref{curvature via torsion and metric})
have the same degree of homogeneity in~$\xi$, namely, $-1$, it is sufficient
to prove formula (\ref{curvature via torsion and metric}) for $\xi$ of norm 1.
Moreover, by rotating our Cartesian coordinate system we can reduce the
case of general $\xi$ of norm 1 to the case
\begin{equation}
\label{special momentum}
\xi=
\begin{pmatrix}
0&0&1
\end{pmatrix}.
\end{equation}

There is one further simplification that can be made: we claim
that it is sufficient to prove
formula (\ref{curvature via torsion and metric}) for the case when
\begin{equation}
\label{special frame}
e_j{}^\alpha(0)=\delta_j{}^\alpha,
\end{equation}
i.e.~for the case when at the point $x=0$
the elements of the frame are aligned with the coordinate axes.
This claim follows from the observation that the
LHS of formula (\ref{curvature via torsion and metric})
is invariant under rigid special unitary transformations of the column-function $v^+$,
$\,v^+\mapsto Rv^+\,$,
where ``rigid'' refers to the fact that the matrix $R\in\mathrm{SU}(2)$ is constant.
Of course, the column-function $Rv^+$ is no longer an eigenvector of the original
principal symbol, but a new principal symbol obtained from the old one by the rigid
special orthogonal transformation of the frame
(\ref{orthogonal transformation of frame 2})
with the $3\times3$ special orthogonal matrix $O$ expressed in terms of the
$2\times2$ special unitary matrix $R$ in accordance with
(\ref{orthogonal transformation of frame 3}).
One can always choose the special unitary matrix $R$ so that
at the point $x=0$
the elements of the new frame are aligned with the coordinate axes
(in fact, there are two possible choices of $R$ which differ by sign).
It remains only to note that direct inspection of formula (\ref{explicit formula for torsion})
shows that torsion is also invariant under rigid special orthogonal transformations of the frame,
and, hence, the tensor $\overset{*}T$ defined by formula (\ref{definition of torsion with a star})
and appearing in the RHS of formula (\ref{curvature via torsion and metric}) is
invariant under rigid special orthogonal transformations of the frame as well.

Having made the simplifying assumptions
(\ref{special metric})--(\ref{special frame}), we are now in a
position to prove formula (\ref{curvature via torsion and metric}).

Let us calculate the RHS of (\ref{curvature via torsion and metric}) first.
In view of (\ref{special frame})
we have, in the linear approximation in $x$,
\begin{equation}
\label{frame in terms of microrotations}
\begin{pmatrix}
e_1{}^1(x)&e_1{}^2(x)&e_1{}^3(x)\\
e_2{}^1(x)&e_2{}^2(x)&e_2{}^3(x)\\
e_3{}^1(x)&e_3{}^2(x)&e_3{}^3(x)
\end{pmatrix}
=
\begin{pmatrix}
1&w^3(x)&-w^2(x)\\
-w^3(x)&1&w^1(x)\\
w^2(x)&-w^1(x)&1
\end{pmatrix},
\end{equation}
where $w$ is some smooth vector-function which vanishes at $x=0$.
Formula (\ref{frame in terms of microrotations}) is the standard formula
for the linearization of an orthogonal matrix about the identity;
see also formula (10.1) in \cite{rotational_elasticity}.
Note that in Cosserat elasticity literature the vector-function $w$ is called the
\emph{vector of microrotations}. Substituting
(\ref{frame in terms of microrotations}) into
(\ref{explicit formula for torsion with a star})
and
(\ref{definition of curl})
we get, at $x=0$,
\begin{equation}
\label{torsion with a star in terms of microrotations}
\overset{*}T_{\alpha\beta}=\partial w_\beta/\partial x^\alpha-\delta_{\alpha\beta}\operatorname{div}w,
\end{equation}
which is formula (10.5) from \cite{rotational_elasticity}.
Here we freely lower and raise tensor indices
using the fact that the metric is Euclidean
(in the Euclidean case (\ref{special metric})
it does not matter whether a tensor index
comes as a subscript or a superscript).
Substituting
(\ref{torsion with a star in terms of microrotations})
and
(\ref{special momentum})
into the RHS of (\ref{curvature via torsion and metric}) we get,
at our point $Q\in T'M$,
\begin{equation}
\label{RHS of curvature via torsion and metric}
\frac12\,
\frac
{\overset{*}T{}^{\alpha\beta}\xi_\alpha\xi_\beta}
{(g^{\mu\nu}\xi_\mu\xi_\nu)^{3/2}}\,
=
-\frac12(\partial w^1/\partial x^1+\partial w^2/\partial x^2)\,.
\end{equation}

Let us now calculate the LHS of (\ref{curvature via torsion and metric}).
The equation for
the eigenvector $v^+(x,\xi)$
of the principal symbol is
\begin{equation}
\label{equation for v plus}
\begin{pmatrix}
e_3{}^\alpha\xi_\alpha-\|\xi\|&(e_1-ie_2)^\alpha\xi_\alpha
\\
(e_1+ie_2)^\alpha\xi_\alpha&-e_3{}^\alpha\xi_\alpha-\|\xi\|
\end{pmatrix}
\begin{pmatrix}
v^+_1\\ v^+_2
\end{pmatrix}=0\,.
\end{equation}
In view of
(\ref{special momentum})
and
(\ref{special frame})
the (normalized) solution of
(\ref{equation for v plus})
at our point $Q\in T'M$ is
$v^+=
\begin{pmatrix}
1\\0
\end{pmatrix}$.
Of course, our $v^+(x,\xi)$
is defined up to the gauge transformation
\begin{equation}
\label{gauge transformation of the eigenvector}
v^+\mapsto e^{i\phi^+}v^+,
\end{equation}
where
\begin{equation}
\label{phase appearing in gauge transformation}
\phi^+:T'M\to\mathbb{R}
\end{equation}
is an arbitrary smooth function,
however the LHS of (\ref{curvature via torsion and metric})
is invariant under this gauge transformation.
We now perturb equation
(\ref{equation for v plus})
about the point $Q\in T'M$,
that is, about
$x=0$,
$\,\xi=
\begin{pmatrix}
0&0&1
\end{pmatrix}$,
making use of formula
(\ref{frame in terms of microrotations}),
which gives us the following equation for the
increment $\delta v^+$ of
the eigenvector $v^+(x,\xi)$
of the principal symbol:
\begin{multline*}
\begin{pmatrix}
0&0
\\
0&-2
\end{pmatrix}
\begin{pmatrix}
\delta v^+_1
\\
\delta v^+_2
\end{pmatrix}
+
\begin{pmatrix}
0&-w^2(x)-iw^1(x)
\\
-w^2(x)+iw^1(x)&0
\end{pmatrix}
\begin{pmatrix}
1
\\
0
\end{pmatrix}
\\
+
\begin{pmatrix}
0&\delta\xi_1-i\delta\xi_2
\\
\delta\xi_1+i\delta\xi_2&-2\delta\xi_3
\end{pmatrix}
\begin{pmatrix}
1
\\
0
\end{pmatrix}
=0,
\end{multline*}
or, equivalently,
\begin{equation}
\label{equation for increment of v plus component 2}
\delta v^+_2
=
\frac12
(-w^2(x)+iw^1(x)+\delta\xi_1+i\delta\xi_2).
\end{equation}
Formula
(\ref{equation for increment of v plus component 2})
has to be supplemented by the normalization condition
\linebreak
$\|v^+(x,\xi)\|=1$, which in its linearized form reads
\begin{equation}
\label{equation for increment of v plus component 1}
\operatorname{Re}
\delta v^+_1=0.
\end{equation}
Formulae
(\ref{equation for increment of v plus component 1})
and
(\ref{equation for increment of v plus component 2})
define $\delta v^+$
modulo an arbitrary
$\operatorname{Im}
\delta v^+_1$, with this degree of freedom being associated
with the gauge transformation
(\ref{gauge transformation of the eigenvector}),
(\ref{phase appearing in gauge transformation}).
Without loss of generality we may assume that the gauge is
chosen so that
\begin{equation}
\label{equation for increment of v plus component 1 choice of gauge}
\operatorname{Im}
\delta v^+_1=0.
\end{equation}

Combining formulae
(\ref{equation for increment of v plus component 1}),
(\ref{equation for increment of v plus component 1 choice of gauge}) and
(\ref{equation for increment of v plus component 2})
we get
\begin{equation}
\label{formula for increment of v plus}
\delta v^+=
\frac12
\begin{pmatrix}
0
\\
-w^2(x)+iw^1(x)+\delta\xi_1+i\delta\xi_2
\end{pmatrix}.
\end{equation}
Recall that the $w$ appearing in this formula
is some smooth vector-function which vanishes at $x=0$.

Differentiation of
(\ref{formula for increment of v plus})
gives us
\begin{equation}
\label{derivative of v plus in x}
\frac{\partial v^+}{\partial x^\alpha}=
\frac12
\begin{pmatrix}
0
\\
-\partial w^2/\partial x^\alpha+i\partial w^1/\partial x^\alpha
\end{pmatrix},
\end{equation}
\begin{equation}
\label{derivative of v plus in xi}
\frac{\partial v^+}{\partial\xi_1}=
\frac12
\begin{pmatrix}
0
\\
1
\end{pmatrix},
\qquad
\frac{\partial v^+}{\partial\xi_2}=
\frac12
\begin{pmatrix}
0
\\
i
\end{pmatrix},
\qquad
\frac{\partial v^+}{\partial\xi_3}=
0.
\end{equation}
Formulae
(\ref{derivative of v plus in x})
and
(\ref{derivative of v plus in xi})
imply that at our point $Q\in T'M$
\begin{equation}
\label{LHS of curvature via torsion and metric}
-i\{[v^+]^*,v^+\}
=
-\frac12(\partial w^1/\partial x^1+\partial w^2/\partial x^2).
\end{equation}

Comparing formulae
(\ref{RHS of curvature via torsion and metric})
and
(\ref{LHS of curvature via torsion and metric})
and recalling (\ref{topological invariant equals one}),
we arrive at the required result
(\ref{curvature via torsion and metric}).
$\square$

\section{Integration of the curvature term}
\label{Integration}

Substituting
(\ref{curvature via torsion and metric})
into
(\ref{formula for b_2(x) version 2})
we get
\begin{equation}
\label{formula for b_2(x) version 3}
b_2(x)=\frac{9\mathbf{c}}4
\int\limits_{h^+(x,\xi)<1}
\frac
{\overset{*}T{}^{\alpha\beta}(x)\,\xi_\alpha\xi_\beta}
{g^{\mu\nu}(x)\,\xi_\mu\xi_\nu}\,
{d{\hskip-1pt\bar{}}\hskip1pt}\xi\,.
\end{equation}
Recall that $h^+(x,\xi)$ is given by formula (\ref{Hamiltonian expressed via metric}).

The tensor $\overset{*}T$ can be decomposed into pure trace and trace-free pieces, i.e.
\begin{equation}
\label{basic decomposition of torsion}
\overset{*}T{}^{\alpha\beta}=\frac13g^{\alpha\beta}\overset{*}T{}^\gamma{}_\gamma
+
\left(
\overset{*}T{}^{\alpha\beta}
-
\frac13g^{\alpha\beta}\overset{*}T{}^\gamma{}_\gamma
\right).
\end{equation}
It is easy to see that the trace-free piece
(second term in the RHS of (\ref{basic decomposition of torsion}))
does not contribute to the integral in (\ref{formula for b_2(x) version 3}),
hence formula (\ref{formula for b_2(x) version 3}) becomes
\begin{equation}
\label{formula for b_2(x) version 4}
b_2(x)=\frac{3\mathbf{c}}4
\,\overset{*}T{}^\gamma{}_\gamma(x)\,
\int\limits_{h^+(x,\xi)<1}
{d{\hskip-1pt\bar{}}\hskip1pt}\xi
\,=\frac{\mathbf{c}}{8\pi^2}
\bigl(\,\overset{*}T{}^\gamma{}_\gamma
\,\sqrt{\det g_{\alpha\beta}}\,\bigr)(x)\,.
\end{equation}

But formulae
(\ref{definition of axial torsion}),
(\ref{definition of Hodge star})
and
(\ref{definition of torsion with a star})
imply that
\begin{equation}
\label{relation between trace of star T and Hodge dual of axial torsion}
\overset{*}T{}^\gamma{}_\gamma=3*T^\mathrm{ax}.
\end{equation}
Combining formulae
(\ref{formula for b(x) sum of two}),
(\ref{formula for b_1(x) version 4}),
(\ref{formula for b_2(x) version 4})
and
(\ref{relation between trace of star T and Hodge dual of axial torsion})
we arrive at formula
(\ref{formula for b(x)}).
This completes the proof of Theorem \ref{theorem 1.1}.

\section{The subprincipal symbol of the massless Dirac operator}
\label{The subprincipal symbol}

In this section we calculate the subprincipal symbol of the massless Dirac operator,
which prepares the ground for the proof of Theorem \ref{theorem 1.2} in the next
section. In view of Remark 2.1.10 from \cite{mybook}, defining the subprincipal
symbol for the massless Dirac operator on spinors (\ref{definition of Weyl operator}) is problematic,
hence, we work with the massless Dirac operator on half-densities
(\ref{definition of Weyl operator on half-densities}).
For the sake of brevity we denote
the massless Dirac operator on half-densities by $A$ rather than by $W_{1/2}\,$.

\begin{lemma}
\label{The subprincipal symbol lemma}
The subprincipal symbol of the
massless Dirac operator on half-densities (\ref{definition of Weyl operator on half-densities}) is
\begin{equation}
\label{Part 1 of the proof of Theorem eq 1}
A_\mathrm{sub}(x)=
\frac{3\mathbf{c}}4
\,\bigl(*T^\mathrm{ax}(x)\bigr)\,I\,,
\end{equation}
where $\mathbf{c}=\pm1$ is the topological charge
(\ref{definition of relative orientation more natural}),
$T^\mathrm{ax}$ is axial torsion (\ref{definition of axial torsion}),
 $*$ is the Hodge star (\ref{definition of Hodge star})
and $I$ is the $2\times2$ identity matrix.
\end{lemma}

\emph{Proof\ }
We give the proof of (\ref{Part 1 of the proof of Theorem eq 1}) for the case
(\ref{topological invariant equals one}).
There is no need to give a separate proof for
the case $\mathbf{c}=-1$ as the two cases reduce to one another
by inversion of the frame:
the full symbol of the massless Dirac operator on half-densities
changes sign under inversion of the frame
and hence its subprincipal symbol
changes sign under inversion of the frame,
whereas torsion
(\ref{explicit formula for torsion})
is invariant under inversion of the frame.

According to formula (1.2) from \cite{DuiGui}
the subprincipal symbol is defined as
\begin{equation}
\label{definition of subprincipal symbol}
A_\mathrm{sub}:=
A_0+\frac i2
(A_1)_{x^\alpha\xi_\alpha}\,,
\end{equation}
where $A_1(x,\xi)$ and $A_0(x)$ are the homogeneous (in $\xi$)
components of the full symbol
$A(x,\xi)=A_1(x,\xi)+A_0(x)$ of our first order differential operator,
with the subscript indicating degree of homogeneity.
For the massless Dirac operator on half-densities (\ref{definition of Weyl operator on half-densities})
these homogeneous components read (\ref{principal symbol via Pauli matrices}) and
\begin{equation}
\label{formula for A0}
A_0(x)=
-\frac i4\sigma^\alpha
\sigma_\beta
\left(
\frac{\partial\sigma^\beta}{\partial x^\alpha}
+\left\{{{\beta}\atop{\alpha\gamma}}\right\}\sigma^\gamma
\right)
+\frac i2\sigma^\alpha
\left\{{{\beta}\atop{\alpha\beta}}\right\}
\end{equation}
respectively.
Note that in writing down (\ref{formula for A0})
we used the standard formula
\[
\frac1{2\det g_{\kappa\lambda}}\frac{\partial\det g_{\mu\nu}}{\partial x^\alpha}
=\left\{{{\beta}\atop{\alpha\beta}}\right\}.
\]
Our task is to substitute
(\ref{principal symbol via Pauli matrices})
and
(\ref{formula for A0})
into (\ref{definition of subprincipal symbol}).

We fix an arbitrary point $P\in M$ and prove formula
(\ref{Part 1 of the proof of Theorem eq 1}) at this point.
As the LHS and RHS of (\ref{Part 1 of the proof of Theorem eq 1}) are invariant
under changes of local coordinates~$x$,
it is sufficient to check the identity (\ref{Part 1 of the proof of Theorem eq 1}) in
Riemann normal coordinates, i.e.~local coordinates
such that $x=0$ corresponds to the point $P$,
$g_{\mu\nu}(0)=\delta_{\mu\nu}$ and $\frac{\partial g_{\mu\nu}}{\partial x^\lambda}(0)=0$.
Moreover, as the principal symbol is linear in $\xi$ and the formula we are proving involves only
first partial derivatives in $x$, we may assume, without loss of generality,
that we have (\ref{special metric}) for all $x$ in some neighbourhood of the origin.
In other words, it is sufficient to prove formula
(\ref{Part 1 of the proof of Theorem eq 1})
for the case of Euclidean metric.
Furthermore, by rotating our Cartesian coordinate system we can
achieve (\ref{special frame}), which opens the way to the use,
in the linear approximation in $x$,
of formula (\ref{frame in terms of microrotations}).

Substituting (\ref{frame in terms of microrotations}) into
(\ref{Pauli matrices 1}), we get,
in the linear approximation in $x$,
\begin{multline}
\label{Part 1 of the proof of Theorem eq 2}
\sigma^1=
\begin{pmatrix}
w^2&1+iw^3\\
1-iw^3&-w^2
\end{pmatrix}
=\sigma_1\,,
\\
\sigma^2=
\begin{pmatrix}
-w^1&-i+w^3\\
i+w^3&w^1
\end{pmatrix}
=\sigma_2\,,
\\
\sigma^3=
\begin{pmatrix}
1&-iw^1-w^2\\
iw^1-w^2&-1
\end{pmatrix}
=\sigma_3\,.
\end{multline}
Recall that the $w$ appearing in this formula
is some smooth vector-function which vanishes at $x=0$.

Substitution of
(\ref{Part 1 of the proof of Theorem eq 2})
into
(\ref{principal symbol via Pauli matrices})
and
(\ref{formula for A0})
gives us
\begin{multline}
\label{Part 1 of the proof of Theorem eq 3}
A_1(x,\xi)=
\begin{pmatrix}
\xi_3&\xi_1-i\xi_2\\
\xi_1+i\xi_2&-\xi_3
\end{pmatrix}
\\
+
\begin{pmatrix}
w^2\xi_1-w^1\xi_2&iw^3\xi_1+w^3\xi_2+(-iw^1-w^2)\xi_3\\
-iw^3\xi_1+w^3\xi_2+(iw^1-w^2)\xi_3&-w^2\xi_1+w^1\xi_2
\end{pmatrix},
\end{multline}
\begin{equation}
\label{Part 1 of the proof of Theorem eq 4}
A_0(0)=-\frac i4
\begin{pmatrix}
0&1\\
1&0
\end{pmatrix}
\begin{pmatrix}
0&1\\
1&0
\end{pmatrix}
\begin{pmatrix}
\partial w^2/\partial x^1&i\partial w^3/\partial x^1\\
-i\partial w^3/\partial x^1&-\partial w^2/\partial x^1
\end{pmatrix}+\ldots.
\end{equation}
Here formula
(\ref{Part 1 of the proof of Theorem eq 3})
is written in the linear approximation in $x$,
whereas formula
(\ref{Part 1 of the proof of Theorem eq 4})
displays, for the sake of brevity, only one term out of nine
(the one corresponding to $\alpha=\beta=1$ in
(\ref{formula for A0})), with the remaining
eight terms concealed within the dots $\ldots$.
Note also that the Christoffel symbols
disappeared because of our assumption
that the metric is Euclidean.

Substituting
(\ref{Part 1 of the proof of Theorem eq 4})
and
(\ref{Part 1 of the proof of Theorem eq 3})
into
(\ref{definition of subprincipal symbol}),
we get
\begin{equation}
\label{Part 1 of the proof of Theorem eq 5}
A_\mathrm{sub}(0)=
-\frac12\,
(\operatorname{div}w)\,I.
\end{equation}
But, according to
(\ref{torsion with a star in terms of microrotations}),
\begin{equation}
\label{Part 1 of the proof of Theorem eq 6}
\overset{*}T{}^\gamma{}_\gamma(0)=-2\operatorname{div}w.
\end{equation}
Formulae
(\ref{Part 1 of the proof of Theorem eq 5}),
(\ref{Part 1 of the proof of Theorem eq 6}),
(\ref{relation between trace of star T and Hodge dual of axial torsion})
and (\ref{topological invariant equals one})
imply formula
(\ref{Part 1 of the proof of Theorem eq 1})
at $x=0$.~$\square$

\section{Proof of Theorem \ref{theorem 1.2}}
\label{Proof of Theorem}

As Theorem \ref{theorem 1.2} is an if and only if theorem,
our proof comes in two parts.

\

\emph{Part 1 of the proof\ }
Let $A$ be a massless Dirac operator on half-densities. We need to
prove that
a)~the subprincipal symbol of this operator,
$A_\mathrm{sub}(x)$, is proportional to the identity matrix and b)
the second asymptotic coefficient of the spectral function, $b(x)$,
is zero. The required result follows from
Lemma~\ref{The subprincipal symbol lemma}
and
Theorem~\ref{theorem 1.1}.

\

\emph{Part 2 of the proof\ }
Let $A$ be a differential operator such that
a) the subprincipal symbol of this operator,
$A_\mathrm{sub}(x)$, is proportional to the identity matrix and b)
the second asymptotic coefficient of the spectral function, $b(x)$,
is zero. We need to prove that $A$ is a massless Dirac operator on half-densities.

Theorem ~\ref{theorem 1.1} implies that the subprincipal symbol of our operator $A$
is given by formula (\ref{Part 1 of the proof of Theorem eq 1}).
Let $e_j$ be the frame corresponding to the principal
symbol of the operator $A$,
see formulae
(\ref{principal symbol via Pauli matrices})
and
(\ref{frame via Pauli matrices}).
Now, let $B$ be the massless Dirac operator on half-densities
corresponding to the same frame. Then the principal symbols of the
operators $A$ and $B$ coincide.
But Lemma~\ref{The subprincipal symbol lemma} implies that
the subprincipal symbols of the
operators $A$ and $B$ coincide as well.
A first order differential operator is determined by its
principal and subprincipal symbols, hence, $A=B$.~$\square$

\section{Explicit formula for axial torsion}
\label{Explicit formula}

Torsion is a rank three tensor antisymmetric in the last two indices.
It is known \cite{rotational_elasticity,cartantorsionreview} that
torsion has three irreducible pieces.
Only one of the three irreducible pieces of torsion, namely,
the piece which theoretical physicists label by the adjective ``axial'',
appears in our spectral theoretic results,
see Theorem~\ref{theorem 1.1} and Lemma~\ref{The subprincipal symbol lemma}.
It is also interesting that axial torsion is the
irreducible piece which is used when one models the
massless neutrino~\cite{MR2670535} or the electron \cite{mathematika} by means of Cosserat elasticity.

Axial torsion is defined as the totally antisymmetric piece of the torsion tensor,
see formula (\ref{definition of axial torsion}).
This means that axial torsion is a 3-form.
In view of the importance of axial torsion,
we give an explicit formula for its Hodge dual in terms of the
principal symbol $A_1(x,\xi)$.
Formulae
(\ref{explicit formula for torsion with a star}),
(\ref{definition of curl})
and
(\ref{relation between trace of star T and Hodge dual of axial torsion})
imply
\begin{multline}
\label{explicit formula for the trace of torsion with a star}
*T^\mathrm{ax}
=\frac{\delta_{kl}}3\,\sqrt{\det g^{\alpha\beta}}
\ \bigl[
e^k{}_{1}
\,\partial e^l{}_{3}/\partial x^2
+
e^k{}_{2}
\,\partial e^l{}_{1}/\partial x^3
+
e^k{}_{3}
\,\partial e^l{}_{2}/\partial x^1
\\
-
e^k{}_{1}
\,\partial e^l{}_{2}/\partial x^3
-
e^k{}_{2}
\,\partial e^l{}_{3}/\partial x^1
-
e^k{}_{3}
\,\partial e^l{}_{1}/\partial x^2
\bigr].
\end{multline}
Here the coframe $e^k$
is determined from the principal symbol
in accordance with formulae
(\ref{principal symbol via Pauli matrices}),
(\ref{frame via Pauli matrices})
and
(\ref{alternative definition of coframe}),
whereas the contravariant metric tensor $g^{\alpha\beta}$
is determined from the principal symbol
in accordance with formula
(\ref{definition of metric}).

\section*{Acknowledgements}

The authors are grateful to D.R.~Heath-Brown for helpful advice.

\appendix

\section{The massless Dirac operator}
\label{The massless Dirac operator}

Let $M$ be a 3-dimensional connected compact oriented manifold equipped with a
Riemannian metric $g_{\alpha\beta}$, $\alpha,\beta=1,2,3$ being the
tensor indices. Note that we are more prescriptive in this appendix
than in the main text of the paper: in the main text orientability
emerged as a consequence of the
existence of a principal symbol
and the metric was defined via the principal symbol,
whereas in this appendix orientability and metric are introduced \emph{a priori}.

We work only in local coordinates with prescribed orientation.

It is known \cite{Stiefel,Kirby} that
a 3-dimensional oriented manifold is parallelizable,
i.e.~there exist smooth real vector fields $e_j$, $j=1,2,3$, that
are linearly independent at every point $x$ of the manifold.
(This fact is often referred to as \emph{Steenrod's theorem}.)
Each vector $e_j(x)$ has coordinate components $e_j{}^\alpha(x)$,
$\alpha=1,2,3$. Note that we use the Latin letter $j$ for
enumerating the vector fields (this is an \emph{anholonomic} or
\emph{frame} index) and the Greek letter $\alpha$ for enumerating
their components (this is a \emph{holonomic} or \emph{tensor}
index). The triple of linearly independent vector fields $e_j$,
$j=1,2,3$, is called a \emph{frame}. Without loss of generality we
assume further on that the vector fields $e_j$ are orthonormal
with respect to our metric: this can always be
achieved by means of the Gram--Schmidt process.

Define Pauli matrices
\begin{equation}
\label{Pauli matrices 1}
\sigma^\alpha(x):=
s^j\,e_j{}^\alpha(x)\,,
\end{equation}
where
\begin{equation}
\label{Pauli matrices 2}
s^1:=
\begin{pmatrix}
0&1\\
1&0
\end{pmatrix}
=s_1
\,,
\quad
s^2:=
\begin{pmatrix}
0&-i\\
i&0
\end{pmatrix}
=s_2
\,,
\quad
s^3:=
\begin{pmatrix}
1&0\\
0&-1
\end{pmatrix}
=s_3
\,.
\end{equation}
In formula (\ref{Pauli matrices 1})
summation is carried out over the repeated frame index $j=1,2,3$,
and $\alpha=1,2,3$ is the free tensor index.

The massless Dirac operator is the matrix operator
\begin{equation}
\label{definition of Weyl operator}
W:=-i\sigma^\alpha
\left(
\frac\partial{\partial x^\alpha}
+\frac14\sigma_\beta
\left(
\frac{\partial\sigma^\beta}{\partial x^\alpha}
+\left\{{{\beta}\atop{\alpha\gamma}}\right\}\sigma^\gamma
\right)
\right),
\end{equation}
where summation is carried out over
$\alpha,\beta,\gamma=1,2,3$, and
\begin{equation}
\label{definition of Christoffel symbols}
\left\{{{\beta}\atop{\alpha\gamma}}\right\}:=
\frac12g^{\beta\delta}
\left(
\frac{\partial g_{\gamma\delta}}{\partial x^\alpha}
+
\frac{\partial g_{\alpha\delta}}{\partial x^\gamma}
-
\frac{\partial g_{\alpha\gamma}}{\partial x^\delta}
\right)
\end{equation}
are the Christoffel symbols.
Here and throughout this appendix we raise and lower
tensor indices using the metric.
Note that we chose the letter ``$W$'' for denoting the massless Dirac operator
because in theoretical physics literature it is often referred to as the \emph{Weyl}
operator.

Formula (\ref{definition of Weyl operator}) is the formula from
\cite{MR2670535}, only written in matrix notation (i.e.~without spinor
indices). Note that in the process of transcribing formulae from
\cite{MR2670535} into matrix notation we used the identity
\begin{equation}
\label{raising spinor indices in Pauli matrices}
\epsilon\sigma^\alpha\epsilon=(\sigma^\alpha)^T,
\end{equation}
$\alpha=1,2,3$, where
\begin{equation}
\label{metric spinor}
\epsilon:=
\begin{pmatrix}
0&-1\\
1&0
\end{pmatrix}
\end{equation}
is the ``metric spinor''.
The identity (\ref{raising spinor indices in Pauli matrices})
gives a simple way of raising/lowering spinor indices in Pauli matrices in the non-relativistic
($\alpha\ne0$) setting.

Our definition (\ref{definition of Weyl operator}) of the massless Dirac operator
is a special case of the definition from \cite{solovej}.
The two definitions coincide when
we work with a $\mathrm{Spin}$ connection as opposed to a $\mathrm{Spin}^c$ connection,
see Propositions 2.14 and 2.15 in \cite{solovej} for details.

Throughout this paper we work in dimension 3.
The definition of the massless Dirac operator acting over a Riemannian
manifold of arbitrary dimension can be found, for example,
in \cite{gilkey_korea,friedrich,gilkey_asymptotic_book}.

Physically, our massless Dirac operator
(\ref{definition of Weyl operator})
describes a single massless neutrino living in a 3-dimensional compact universe
$M$. The eigenvalues of the massless Dirac operator are the energy levels.


The massless Dirac operator (\ref{definition of Weyl operator}) acts on 2-columns
$v=\begin{pmatrix}v_1&v_2\end{pmatrix}^T$
of complex-valued scalar functions.
In differential geometry this object is referred to as
a (Weyl) spinor so as to emphasize the fact that $v$ transforms
in a particular way under transformations of the orthonormal frame $e_j$.
However, as in our exposition the frame $e_j$ is assumed to be chosen
\emph{a priori}, we can treat the components of the spinor as scalars.
This issue will be revisited below when we state Property 4
of the massless Dirac operator.

We now list the main properties of the massless Dirac operator.
We state these without proofs. The proofs can be found in Appendix 3.A of \cite{ChervovaPhD}
or in \cite{solovej}.

\

\textbf{Property 1.}
The massless Dirac operator is invariant under changes of local
coordinates $x$, i.e.~it maps 2-columns of smooth scalar functions $M\to\mathbb{C}^2$
to 2-columns of smooth scalar functions $M\to\mathbb{C}^2$
regardless of the choice of local coordinates.

\

\textbf{Property 2.}
The massless Dirac operator
is formally self-adjoint (symmetric)
with respect to the inner product
\begin{equation}
\label{inner product on colums of scalars}
\int_Mw^*v\,\sqrt{\det g_{\alpha\beta}}\,dx
\end{equation}
on 2-columns of smooth scalar functions $v,w:M\to\mathbb{C}^2$.

\

\textbf{Property 3.}
The massless Dirac operator $W$ commutes
\begin{equation}
\label{commutes}
\mathrm{C}(Wv)=W\mathrm{C}(v)
\end{equation}
with the antilinear map
\begin{equation}
\label{antilinear}
v\mapsto\mathrm{C}(v):=\epsilon\overline{v},
\end{equation}
where $\epsilon$ is the ``metric spinor'' (\ref{metric spinor}).
In theoretical physics the transformation
(\ref{antilinear})
is referred to as \emph{charge conjugation} \cite{PCT,solovej}.

\

Formula (\ref{commutes}) implies that $v$ is an eigenfunction
of the massless Dirac operator corresponding to an eigenvalue $\lambda$
if and only if
$\mathrm{C}(v)$ is an eigenfunction
of the massless Dirac operator corresponding to the same eigenvalue
$\lambda$.
Hence, all eigenvalues of the massless Dirac operator have even
multiplicity. Moreover, any eigenfunction $v$ and its ``partner'' $\mathrm{C}(v)$
make the same contribution to the spectral function
(\ref{definition of spectral function})
at every point $x$ of the manifold $M$.

If, as in \cite{solovej}, we introduce a magnetic field, then we lose the commutation property
(\ref{commutes}) and the double eigenvalues split up. This indicates that the double eigenvalues
of the massless Dirac operator correspond to the two different spins.

\

\textbf{Property 4.}
This property has to do with a particular behaviour under $\mathrm{SU}(2)$ transformations.
Let $R:M\to\mathrm{SU}(2)$ be an arbitrary smooth special unitary matrix-function.
Let us introduce new Pauli matrices
\begin{equation}
\label{special unitary transformation of Pauli matrices}
\tilde\sigma^\alpha:=R\sigma^\alpha R^*
\end{equation}
and a new operator $\tilde W$ obtained by replacing the $\sigma$
in (\ref{definition of Weyl operator}) by $\tilde\sigma$.
It turns out (and this is Property 4) that the two operators,
$\tilde W$ and $W$, are related in exactly the same way as the Pauli
matrices, $\tilde\sigma$ and $\sigma$, that is,
\begin{equation}
\label{special unitary transformation of Weyl operator}
\tilde W=RW R^*.
\end{equation}

\

We now examine the geometric meaning of the
transformation (\ref{special unitary transformation of Pauli matrices}).
Let us expand the new Pauli matrices $\tilde\sigma$ with respect to
the basis (\ref{Pauli matrices 2}):
\begin{equation}
\label{Pauli matrices 1 new}
\tilde\sigma^\alpha(x)=
s^j\,\tilde e_j{}^\alpha(x).
\end{equation}
Formulae
(\ref{Pauli matrices 1}),
(\ref{Pauli matrices 1 new})
and
(\ref{special unitary transformation of Pauli matrices})
give us the following identity relating
the new vector fields $\tilde e_j$
and the old vector fields $e_j$:
\begin{equation}
\label{orthogonal transformation of frame 1}
Rs^kR^*e_k=s^j\tilde e_j\,.
\end{equation}
Resolving (\ref{orthogonal transformation of frame 1})
for $\tilde e_j$ we get
\begin{equation}
\label{orthogonal transformation of frame 2}
\tilde e_j=O_j{}^ke_k\,,
\end{equation}
where the real scalars $O_j{}^k$ are given by the formula
\begin{equation}
\label{orthogonal transformation of frame 3}
O_j{}^k=\frac12\operatorname{tr}(s_jRs^kR^*)\,.
\end{equation}
Note that in writing formulae
(\ref{orthogonal transformation of frame 1}) and (\ref{orthogonal transformation of frame 2})
we chose to hide the tensor index, i.e.~we chose to hide the coordinate components of our vector fields.
Say, formula (\ref{orthogonal transformation of frame 2}) written
in more detailed form reads $\tilde e_j{}^\alpha=O_j{}^ke_k{}^\alpha$.

The scalars (\ref{orthogonal transformation of frame 3})
can be viewed as elements of a real $3\times3$ matrix-function $O$
with the first index, $j$, enumerating rows and the second, $k$, enumerating columns.
It is easy to check that this matrix-function $O$ is special orthogonal.
Hence, the new vector fields $\tilde e_j$ are orthonormal and have the same orientation
as the old vector fields $e_j$.
We have shown that the
transformation (\ref{special unitary transformation of Pauli matrices})
has the geometric meaning of switching from our original oriented orthonormal frame $e_j$
to a new oriented orthonormal frame $\tilde e_j$.

Formula (\ref{orthogonal transformation of frame 3}) means that the special
unitary matrix $R$ is, effectively, a square root of the special orthogonal
matrix $O$.
It is easy to see that for a given matrix $O\in\mathrm{SO}(3)$ formula
(\ref{orthogonal transformation of frame 3}) defines the
matrix $R\in\mathrm{SU}(2)$ uniquely up to sign.
This observation allows us to view the issue
of the geometric meaning of the
transformation (\ref{special unitary transformation of Pauli matrices})
the other way round:
given a pair of orthonormal frames, $e_j$ and $\tilde e_j$, with the
same orientation
(i.e.~with $\operatorname{sgn}\det e_j{}^\alpha=\operatorname{sgn}\det\tilde e_j{}^\alpha$),
we can recover the special orthogonal
matrix-function $O(x)$ from formula (\ref{orthogonal transformation of frame 2})
and then attempt finding a smooth special unitary matrix-function $R(x)$
satisfying (\ref{orthogonal transformation of frame 3}).
Unfortunately, this may not always be possible due to topological obstructions.
We can only guarantee the absence of topological obstructions when the two
frames, $e_j$ and $\tilde e_j$, are sufficiently close to each other, which
is equivalent to saying that
we can only guarantee the absence of topological obstructions when
the special orthogonal matrix-function $O(x)$
is sufficiently close to the identity matrix for all $x\in M$.

We illustrate the possibility of a topological obstruction by means of an explicit example.
Consider the unit torus $\mathbb{T}^3$ parameterized by cyclic coordinates $x^\alpha$,
$\alpha=1,2,3$, of period $2\pi$. The metric is assumed to be Euclidean. Define a pair
of orthonormal frames
\begin{equation}
\label{frame on torus standard}
e_j{}^\alpha:=\delta_j{}^\alpha
\end{equation}
and
\begin{equation}
\label{frame on torus}
\tilde e_1{}^\alpha:=
\begin{pmatrix}
\cos k_3x^3\\
\sin k_3x^3\\
0
\end{pmatrix},
\qquad
\tilde e_2{}^\alpha:=
\begin{pmatrix}
-\sin k_3x^3\\
\cos k_3x^3\\
0
\end{pmatrix},
\qquad
\tilde e_3{}^\alpha:=
\begin{pmatrix}
0\\
0\\
1
\end{pmatrix},
\end{equation}
where $k_3$ is an odd integer.
Let $W$ and $\tilde W$ be the massless Dirac operators corresponding to the frames
(\ref{frame on torus standard})
and
(\ref{frame on torus}) respectively.
We claim that there does not exist
a smooth matrix-function $R:\mathbb{T}^3\to\mathrm{SU}(2)$ which would give
(\ref{orthogonal transformation of frame 3}),
where $O(x)$ is the special orthogonal matrix-function defined by formula
(\ref{orthogonal transformation of frame 2}).
We justify this claim in two different ways.

\emph{Justification 1.}
Resolving the system (\ref{orthogonal transformation of frame 2})--(\ref{frame on torus})
locally for $R$, we get
\begin{equation}
\label{example of topological obstruction formula for R}
R(x^3)=
\pm
\begin{pmatrix}
e^{\frac i2k_3x^3}&0\\
0&e^{-\frac i2k_3x^3}
\end{pmatrix},
\end{equation}
and this solution is unique modulo choice of sign;
here the freedom in the choice of sign is not surprising
as $\mathrm{SU}(2)$ is the double cover of $\mathrm{SO}(3)$.
Formula (\ref{example of topological obstruction formula for R})
defines a continuous single-valued matrix-function on the unit torus $\mathbb{T}^3$
if and only if the integer $k_3$ is even, which it is not.

\emph{Justification 2.}
It is sufficient to show that the two operators,
$W$ and $\tilde W$, have different spectra.
Straightforward separation of variables shows that zero is an eigenvalue of the operator
$W$ but not an eigenvalue of the operator $\tilde W$.

One can generalize the above example
by introducing rotations in three different directions, which leads
to eight genuinely distinct parallelizations. See also \cite{wild} page~524 or \cite{bar_2000} page~21.

Let us emphasize that the topological obstructions we were discussing
have nothing to do with Stiefel--Whitney classes. We are working on
a parallelizable manifold and the Stiefel--Whitney class of such
a manifold is trivial. The topological issue at hand is that our
parallelizable manifold
may be equipped with different spin structures.

We say that two massless Dirac operators, $W$ and $\tilde W$, are
equivalent if there exists a smooth matrix-function
$R:M\to\mathrm{SU}(2)$ such that the corresponding Pauli matrices,
$\sigma^\alpha$ and $\tilde\sigma^\alpha$, are related in accordance
with (\ref{special unitary transformation of Pauli matrices}).
In view of Property 4 (see formula
(\ref{special unitary transformation of Weyl operator}))
all massless Dirac operators from the same
equivalence class generate the same spectral function
(\ref{definition of spectral function})
and the same counting function
(\ref{definition of counting function}),
so for the purposes of our paper viewing such operators as
equivalent is most natural.

As explained above, there may be many distinct equivalence classes
of massless Dirac operators, the difference between
which is topological. Studying the spectral theoretic implications of
these topological differences is beyond the scope of our paper.
The two-term asymptotics
(\ref{two-term asymptotic formula for spectral function for dirac})
and
(\ref{two-term asymptotic formula for counting function for dirac})
derived in the main text of our paper do not feel this topology.

In theoretical physics the $\mathrm{SU}(2)$ freedom involved in defining the
massless Dirac operator is interpreted as a gauge degree of
freedom. We do not adopt this point of view (at least explicitly)
in order to fit the massless Dirac operator into the
standard spectral theoretic framework.

We defined the massless Dirac operator
(\ref{definition of Weyl operator})
as an operator acting on 2-columns of scalar functions,
i.e.~on 2-columns of quantities which do not change under changes of
local coordinates. This necessitated the introduction of the density
$\sqrt{\det g_{\alpha\beta}}$ in the formula
(\ref{inner product on colums of scalars}) for the inner product.
In spectral theory it is more common to work with half-densities.
Hence, we introduce the operator
\begin{equation}
\label{definition of Weyl operator on half-densities}
W_{1/2}:=
(\det g_{\kappa\lambda})^{1/4}
\,W\,
(\det g_{\mu\nu})^{-1/4}
\end{equation}
which maps half-densities to half-densities.
We call the operator
(\ref{definition of Weyl operator on half-densities})
\emph{the massless Dirac operator on half-densities}.

\section{The spectrum for the torus and the sphere}
\label{The torus and the sphere}

In this appendix we examine the massless Dirac operator on the unit torus $\mathbb{T}^3$
and the unit sphere $\mathbb{S}^3$
and compare our asymptotic formulae
(\ref{two-term asymptotic formula for spectral function for dirac})
and
(\ref{two-term asymptotic formula for counting function for dirac})
with known explicit formulae.
The  torus is assumed to be equipped with Euclidean metric
(see also Appendix \ref{The massless Dirac operator})
whereas the sphere is assumed to be equipped with metric induced by the natural embedding
of $\mathbb{S}^3$ in Euclidean space $\mathbb{R}^4$.
Note that in view of the obvious symmetries of the torus and the sphere the scalar function
$e(\lambda,x,x)/\sqrt{\det g_{\alpha\beta}(x)}$ is constant
(see also Remark~\ref{remark on spectral function desnity versus scalar}),
so formulae
(\ref{two-term asymptotic formula for spectral function for dirac})
and
(\ref{two-term asymptotic formula for counting function for dirac})
are in this case equivalent, in the sense that they follow from one another.
Hence, we will be dealing with formula
(\ref{two-term asymptotic formula for counting function for dirac}) only.

We have $\operatorname{Vol}\mathbb{T}^3=(2\pi)^3$, so for the torus formula
(\ref{two-term asymptotic formula for counting function for dirac})
reads
\begin{equation}
\label{two-term asymptotic formula for counting function for dirac torus}
N(\lambda)=\frac43\pi\lambda^3+o(\lambda^2).
\end{equation}
The nonperiodicity condition (see Definitions 8.3 and 8.4 in \cite{jst_part_a})
is fulfilled for the torus, so,
according to Theorem 8.4 from \cite{jst_part_a},
the asymptotic formula
(\ref{two-term asymptotic formula for counting function for dirac torus})
holds as it is, without mollification.

In order to test formula
(\ref{two-term asymptotic formula for counting function for dirac torus})
we calculate the spectrum of the massless Dirac operator on
$\mathbb{T}^3$ explicitly.
We do this first for the spin structure associated with the frame
(\ref{frame on torus standard}). Then the spectrum is as follows.
\begin{itemize}
\item
Zero is an eigenvalue of multiplicity two.
\item
For each $m\in\mathbb{Z}^3\setminus\{0\}$ we have
the eigenvalue $\|m\|$ and unique (up to rescaling) eigenfunction,
with eigenfunctions corresponding to different $m$ being linearly
independent.
\item
For each $m\in\mathbb{Z}^3\setminus\{0\}$ we have
the eigenvalue $\,-\|m\|\,$ and unique (up to rescaling) eigenfunction,
with eigenfunctions corresponding to different $m$ being linearly
independent.
\end{itemize}
Hence, $N(\lambda)+1$ is the number of integer lattice points
inside a 2-sphere of radius $\lambda$ in $\mathbb{R}^3$ centred at the origin.
According to \cite{heath-brown_1999} the latter
admits the asymptotic expansion
\begin{equation}
\label{heath-browns_result}
\frac43\pi\lambda^3+O_\varepsilon(\lambda^{21/16+\varepsilon})
\end{equation}
as $\lambda\to+\infty$, with $\varepsilon$ being an arbitrary positive number.
This agrees with our asymptotic formula
(\ref{two-term asymptotic formula for counting function for dirac torus}).

Let us now consider the
spin structure associated with the frame
(\ref{frame on torus}). Then the spectrum is as follows.
\begin{itemize}
\item
For each $m\in\mathbb{Z}^3$ we have
the eigenvalue $\|m-(0,0,1/2)\|$ and unique (up to rescaling) eigenfunction,
with eigenfunctions corresponding to different $m$ being linearly
independent.
\item
For each $m\in\mathbb{Z}^3$ we have
the eigenvalue $\,-\|m-(0,0,1/2)\|\,$ and unique (up to rescaling) eigenfunction,
with eigenfunctions corresponding to different $m$ being linearly
independent.
\end{itemize}
Hence, $N(\lambda)$ is the number of integer lattice points
inside a 2-sphere of radius $\lambda$ in $\mathbb{R}^3$ centred at $(0,0,1/2)$.
Here the sphere is shifted from the origin so one cannot apply the result from
\cite{heath-brown_1999}.
However, as the shift is rational, one can reduce the problem
to counting integer lattice points in a rational ellipsoid centred at the origin,
and an application of the result from \cite{rational_ellipsoids} gives us for the shifted sphere
the same asymptotic expansion (\ref{heath-browns_result}) as for the sphere centred at the origin.

As explained in Appendix \ref{The massless Dirac operator}, the unit torus $\mathbb{T}^3$
admits a total of eight different spin structures. For each of these the problem of counting
positive eigenvalues of the massless Dirac operator reduces to counting
integer lattice points inside a 2-sphere of radius $\lambda$ in $\mathbb{R}^3$
(possibly, shifted from the origin by a rational shift),
so in all eight cases we do get
(\ref{two-term asymptotic formula for counting function for dirac torus}).
In fact, we can replace
the remainder $o(\lambda^2)$ in (\ref{two-term asymptotic formula for counting function for dirac torus})
by $O_\varepsilon(\lambda^{21/16+\varepsilon})$ and this
holds for all eight different spin structures.

In the remainder of this appendix we examine the massless Dirac operator
on the unit sphere $\mathbb{S}^3$.
We have $\operatorname{Vol}\mathbb{S}^3=2\pi^2$, so for the sphere formula
(\ref{two-term asymptotic formula for counting function for dirac})
reads
\begin{equation}
\label{two-term asymptotic formula for counting function for dirac sphere}
N(\lambda)=\frac{\lambda^3}3+o(\lambda^2).
\end{equation}
The nonperiodicity condition fails for the sphere because all geodesics are closed
with period $2\pi$, so formula
(\ref{two-term asymptotic formula for counting function for dirac sphere})
cannot be used in its original form and has to be mollified,
see Remark \ref{remark on mollification}.
We will deal with the mollification issue later and give explicit formulae for
the eigenvalues first.

It is known that $\mathbb{S}^3$ admits a unique spin structure, see Section 5 in \cite{bar_2000}.
The spectrum of the massless Dirac operator on
$\mathbb{S}^3$ has been computed by different authors using different methods
\cite{sulanke,trautman,bar_1996,bar_2000} and reads as follows:
the eigenvalues are
\begin{equation}
\label{11 August 2012 equation 3}
\pm\left(k+\frac12\right),
\qquad k=1,2,\ldots,
\end{equation}
and their multiplicity is
\begin{equation}
\label{11 August 2012 equation 4}
k(k+1).
\end{equation}

The mollification procedure from Section 7 of \cite{jst_part_a} goes as follows.
Put $N(\lambda):=0$ for $\lambda\le0$ and take an arbitrary
real-valued even function $\rho(\lambda)$ from Schwartz space
$\mathcal{S}(\mathbb{R})$ whose Fourier transform $\hat\rho(t)$
satisfies conditions $\hat\rho(0)=1$ and $\operatorname{supp}\hat\rho\subset(-2\pi,2\pi)$.
Then, according to Theorem 7.2 from \cite{jst_part_a},
the mollified version of formula (\ref{two-term asymptotic formula for counting function for dirac sphere}) reads
\[
\int N(\lambda-\mu)\,\rho(\mu)\,d\mu=\frac{\lambda^3}3+O(\lambda)
\]
and this result holds notwithstanding the failure of the nonperiodicity condition.
However, for the sphere there is a much simpler way of testing our asymptotic formula.
Let $\lambda\ge2$ be integer. Taking an integer $\lambda$ puts us exactly in the middle
of the gap between two consecutive clusters of eigenvalues,
see formulae (\ref{11 August 2012 equation 3}) and (\ref{11 August 2012 equation 4}),
and achieves the same averaging effect as convolution with a function from Schwartz space.
For integer $\lambda\ge2$ we get
\[
N(\lambda)=\sum_{k=1}^{\lambda-1}k(k+1)
=\frac{\lambda^3}3-\frac{\lambda}3
\]
which agrees with our asymptotic formula
(\ref{two-term asymptotic formula for counting function for dirac sphere}).

\end{document}